\documentclass{article}
 \usepackage{ amsthm, amsmath, amssymb}
\usepackage{color}

\theoremstyle{plain}
\newtheorem{thm}{Theorem}[section]
\newtheorem{lem}[thm]{Lemma}
\newtheorem{cor}[thm]{Corollary}
\newtheorem*{centreconj}{Centre Conjecture}

\newtheorem*{partitionprob}{Partition Problem}
\theoremstyle{remark}
\newtheorem{remark}{Remark}
\newtheorem{example}{Example}

\title{The Geometry of an Equifacetal Simplex}
\author{Allan L. Edmonds}
\date{}

\begin{document}

\maketitle
\begin{abstract} Equifacetal simplices, all of whose codimension one faces are congruent to one another, are studied.   It is shown that the isometry group of such a simplex acts transitively on its set of vertices, and, as an application, equifacetal simplices are shown to have unique centres.  It is conjectured that a simplex with a unique centre must be equifacetal. The notion of the combinatorial type of an equifacetal simplex is introduced and analysed, and all possible combinatorial types of equifacetal simplices are constructed in even dimensions. 
\end{abstract}
\section{Introduction}
In this paper we study an appropriate higher-dimensional analogue  of the simplest and most familiar construction of traditional Euclidean plane geometry:
the construction  an equilateral triangle with  any given line segment  as an edge.

We analyse simplices all of whose facets, the maximal proper faces, are congruent to one another, which we call \emph{equifacetal} simplices.  

Moving up a dimension from the classical situation, one may ask, for example, whether for a given triangle there is a tetrahedron all of whose faces are congruent to the given triangle.   It turns out that for a given triangle $T$ there is a tetrahedron $S$ such that all four faces of $S$ are congruent to $T$ if and only if $T$ is acute. Altshiller-Court  \cite{Altshiller1935}, among others, studies such equifacetal tetrahedra.  Later printings of his monograph include the stated restriction on the facets. See also \cite{Gardiner1990, Honsberger1976, HajjaWalker2001}.

In higher dimensions  we show that there are many equifacetal, non-equilateral simplices, systematically constructing examples in a large array of possible ``combinatorial types''.  

In a recent paper, Martini and Wenzel \cite{MartiniWenzel2003} gave a construction that showed the existence of non-regular equifacetal $n$-simplices in higher dimensions.  Devid\'{e}
\cite{Devide1975} presented an analysis that in particular deals with the dimension $4$ case. See also McMullen \cite{McMullen2000} and Weissbach \cite{Weissbach2000} for related work, focused on the requirement that the facets have equal areas.

We study the isometry group of an equifacetal simplex and show that it acts transitively on the vertices and on the facets. As a corollary it follows that an equifacetal simplex has a unique centre.  This result provides a unified explanation of numerous scattered and more or less ad hoc results in the literature where it is shown that certain centres coincide for certain equifacetal simplices (e.g.  \cite{Altshiller1935,Devide1975,Devide1984, HajjaWalker2001, KupitzMartini1994, MartiniWenzel2003})
.

We then consider the converse and conjecture that a simplex with unique centre is equifacetal. In the case of the tetrahedron this follows from work of Altshiller-Court \cite{Altshiller1935}.  See also Devid\'{e} \cite{Devide1975} for proofs that put the result in a more general context. The proof only requires knowing that two of the incentre, circumcentre, centroid or Monge point coincide.  
A proof for the full conjecture should involve systematic use of larger numbers of canonical centres that in general ought to be suitably independent of one another.

\paragraph{Outline of the paper.} 
In Section \ref{sec:tetrahedron} we show that the isometry group of an equifacetal tetrahedron acts transitively on vertices and that an equifacetal tetrahedron has a unique centre. In Section \ref{sec:simplex} we prove the analogous results for equifacetal simplices in general, with somewhat more elaborate proofs required. In Section \ref{sec:combinatorial} we discuss the realization of a given set of edge lengths by a simplex and introduce the strong combinatorial type of a simplex as the edge colouring corresponding to edges of different lengths. In Section \ref{sec:restrictions} we explain the main restrictions on the strong combinatorial types that arise from an equifacetal simplex. In Section \ref{sec:partition} we associate with each equifacetal  $n$-simplex  a partition of $n$, which we call the bookkeeping partition or the weak combinatorial type, as a way of cataloguing different equifacetal simplices. In Section \ref{sec:realizing} we prove certain realization results about which weak combinatorial types arise from strong combinatorial types of equifacetal simplices. In Section \ref{sec:group} we indicate how any transitive permutation group gives rise to an equifacetal simplex. In Section \ref{sec:examples} we display a number of examples in dimensions $5$ and $6$ that temper one's optimism about natural possible extensions of some of the results in low dimensions.  In a final section we list all the realizable weak and strong combinatorial types up through dimension $6$.
\paragraph{Acknowledgement.} 
I would like to thank Horst Martini and Mowaffaq Hajja for their interest in this work and their guidance to the literature.

\section{The Tetrahedron}\label{sec:tetrahedron}

As mentioned above,  Altshiller-Court proved that a triangle is a facet of an equifacetal tetrahedron if and only if it is acute.
Another way to say this is that, given three edge lengths $0<a\le b\le c$, there is a tetrahedron with all faces having these three edge lengths if and only if $c^2<a^2+b^2$, as follows immediately from the Law of Cosines. Such a tetrahedron is characterized by the property  that disjoint edges are congruent to one another.

Altshiller-Court \cite{Altshiller1935} also noted by direct arguments that for an equifacetal tetrahedron (``isosceles'' in his terminology), the centroid, circumcentre, incentre, and Monge point coincide.  He also showed that if any two of these four centres coincide, then the tetrahedron must be equifacetal. Devid\'{e} \cite{Devide1975}, perhaps unaware of Altshiller-Court's work, showed that the incentre, circumcentre, and centroid all coincide for an equifacetal tetrahedron, and that if any two of them coincide, then the tetrahedron is equifacetal. See also \cite{Devide1984} for more aspects of equifacetal tetrahedra. And Kupitz  and Martini \cite{KupitzMartini1994} showed that the Fermat-Torricelli point coincides with the circumcentre  or with the centroid  if and only if the tetrahedron is equifacetal. Hajja and Walker \cite{HajjaWalker2001} also showed that the Fermat-Torricelli point must coincide with the classical centres for equifacetal tetrahedra, and that if the Fermat-Torricelli point coincides with a classical centre, including the incentre, then the tetrahedron must be equifacetal. We will offer an explanation for these and similar results based on the presence of enough symmetry, showing that an equifacetal tetrahedron has just  one unique centre. The separate proof we give here foreshadows a more general result to be presented in due course for arbitrary equifacetal simplices.

Any tetrahedral centre must be fixed by any isometry of the tetrahedron to itself.

\begin{thm}
An equifacetal tetrahedron $T=ABCD$ has \emph{exactly} one centre, for there is a unique fixed point for the group of isometries leaving $T$ invariant.
\end{thm}

\emph{Proof. } 
The equifacetal tetrahedron $T$ has the property that disjoint edges are congruent to one another. Consider the line $L$ connecting the midpoints $M$ of $AB$ and $N$
of
$CD$.  We claim that $L$ is actually perpendicular to $AB$ and $CD$. To see
this, note that $MN$ is the median line of the triangle $ANB$.  But $AN$ is the
median of triangle $ACD$ while $BN$ is the median line of triangle $BDC$ and
these two triangles $ACD$ and $BDC$ are congruent by assumption.  It follows that $AN\cong
BN$.  Therefore triangle $ANB$ is isosceles, and the median line $MN$ of triangle $ANB$ is perpendicular to the
base.

It follows that $T$ is invariant under the $180^\circ$ rotation about line $L$.
And so any tetrahedral centre must lie on $L$. In similar fashion $T$ is
invariant under the $180^\circ$ rotation about the other two bimedian lines. 
These lines have exactly one point in common, the centroid.
\qed\bigskip

The discussion above shows that the converse of the theorem also holds: a tetrahedron with a unique centre must be equifacetal.
\section{The $n$-Simplex}\label{sec:simplex}
By an $n$-simplex we understand the convex hull of $n+1$ affinely independent points $A_{0},A_{1},\dots ,A_{n}$ in some euclidean space $\mathbb{R}^d$, $d\ge n$.  We will write it simply as $A_{0}A_{1}\dots A_{n}$. A point $P$ in the $n$-simplex may be given uniquely by its barycentric coordinates as $P=t_{0}A_{0}+t_{1}A_{1}+\dots t_{n}A_{n}$,
where each $t_{i}\ge 0$ and $\sum t_{i}=1$.

For more details on simplices and polytopes in general we refer to the classic monograph of Gr\"{u}nbaum \cite{Grunbaum1967} or the more recent text by Ziegler \cite{Ziegler1995}.
We shall see that equifacetal simplices in higher dimensions have many of the same properties of their 3-dimensional counterparts.  In particular, their isometry groups act transitively on vertices, so that they have unique centres.  Two differences are that being a face of an equifacetal simplex is no longer an open condition and that many more patterns of equifacetal face types emerge. 
\subsection{Isometry groups}
By the higher dimensional analogue  of the Side-Side-Side principle of Euclidean plane geometry, the edge lengths completely determine a simplex up to isometry, as we formulate in the following lemma, whose formal proof we omit. See Berger \cite{Berger1987}, Proposition 9.7.1, for example, for more details.
\begin{lem}\label{isometrycondition}
Let an $n$-simplex $S=A_0A_1\dots A_n$ be given. If $G$ is a subgroup of  the symmetric group on $n+1$ letters,  acting on the \emph{set} of vertices, then $G$ acts on $S$ by isometries if and only if $G$ preserves distances on
$\{ A_0,A_1,\dots, A_n \}$, where, in barycentric coordinates, we define $g(\sum t_{i}A_{i})=\sum t_{i}g(A_{i})$.
\qed
\end{lem}
\begin{lem}[Vertex Uniformity]\label{graph}
If $S=A_0A_1\dots A_n$ is an equifacetal $n$-simplex, and $\ell$ is an edge length in $S$, then the graph $G_{\ell}$  of edges of  length $\ell$ contains all $n+1$ vertices and all vertices have the same degree.  \end{lem}

\emph{Proof.}
If two vertices had different degrees, then the two facets of $S$ obtained by omitting one or the other of the two vertices would have different numbers of edges of the given length, contradicting the fact that $S$ is equifacetal.  \qed\bigskip

Similar reasoning shows that, in addition, if the graph is not connected, all components of $G_{\ell}$ are congruent to one another.
For,  if two components of $G_{\ell}$ were not congruent, then the two facets obtained by omitting a vertex from the respective components would have non-congruent graphs of edges of the given length. 

The following strengthening of these observations is the most basic and useful result of  the paper.
\begin{thm}[Vertex Transitivity]\label{transitivity}
The isometry group of an equifacetal $n$-simplex  $S=A_0A_1\dots A_n$ acts transitively on the vertices and, equivalently, on the facets.
\end{thm}
\emph{Proof.  } 
The result is clearly true for any equilateral simplex. So assume without loss of generality that the simplex $S$ is not equilateral.

Let $A_i$ and $A_j$ be two vertices.  We seek to display an isometry of $S$ mapping $A_i\to A_j$. Let $F_i$ and $F_j$ be the faces obtained by deleting $A_i$ and $A_j$, respectively.  Let $\varphi:F_i\to F_j$ be any isometry, which exists by assumption.  

We claim that by mapping $A_{i}\to A_{j}$, $\varphi$ extends to an isometry of $S$ to itself.  For each vertex $A_k$ of $F_{i} $ consider the vertex $A_k'=\varphi(A_k)$.  By Lemma \ref{isometrycondition} we just need to verify that $|A_jA_k'|=|A_iA_k|$.  

Since $A_{k}$ and $A_{k}'$ have the same sets (with multiplicities) of edge lengths incident at them in  $F_{i}$ and $F_{j}$, respectively, and the same edge lengths incident at them in $S$ itself, by Lemma \ref{graph}, it follows that the length of the edge $A_{k}A_{i}$ must be the same as the length of the edge $A_{k}'A_{j}$.

Thus, consideration of the constant edge length graphs containing $A_iA_k$, which have constant degree at all vertices of  $S$, implies that $|A_jA_k'|=|A_iA_k|$. 
It follows that the mapping $A_i\to A_j$ extends $\varphi$ to an isometry of all of $S$, as required.
\qed\bigskip

\begin{remark}
Certainly the isometry group of a  non-equilateral simplex is not the full symmetric group.  All of our actual constructions of equifacetal $n$-simplices involve an isometry group that contains a cyclic or dihedral group of order $n+1$. There is a close connection between transitive subgroups of the symmetric group and combinatorial types of equifacetal simplices.  But note that a doubly transitive group of isometries would imply that the simplex is equilateral.
\end{remark}
\subsection{Simplex Centres}
For our purposes it suffices to view a simplex centre as a function that assigns to any simplex a definite point in the affine span of the simplex such that that point is fixed by any isometry of space that leaves the simplex invariant. (If a centre is viewed as a function of the ordered vertices, then one would want it to be invariant under permutation of the vertex labelling. One might also require that it scale appropriately under a similarity map. See Kimberling \cite{Kimberling1997} for such a discussion in dimension $2$.)

For an arbitrary generic simplex literally any point could potentially be a centre.  Since an isometry preserving a simplex must also preserve any centre, the presence of nontrivial symmetry will  restrict the centres of a simplex.

\begin{cor}
An equifacetal $n$-simplex  $S=A_0A_1\dots A_n$ has a unique simplex centre.
\end{cor}
\emph{Proof.  } 
Any isometry must fix the centroid $\sum \frac{1}{n+1}A_i$.  We will show that any point fixed by the entire isometry group of $S$ must be the centroid.

Suppose $g:S\to S$ is an isometry. Then $g$ maps vertices to vertices and preserves barycentric coordinates, so that $g(\sum t_{i}A_{i})=\sum t_{i}g(A_{i})$, where each $t_{i}\ge 0$ and $\sum t_{i}=1$. (From this it follows that $g$ preserves the centroid.) 

Suppose in addition that $x=\sum t_{i}A_{i}$ is a fixed point of $g$. Then we conclude that $\sum t_{i}g(A_{i})=\sum t_{i}A_{i}$. Uniqueness of barycentric coordinates  shows that the coefficient of $A_{i}$ is the same as the coefficient of $g(A_{i})$. If $x$ is a fixed point of the entire isometry group, then $t_{0}=t_{1}=\cdots =t_{n}$ since the isometry group acts transitively on the set of vertices. It follows that $x$ is the centroid $\sum \frac{1}{n+1}A_{i}$.
\qed\bigskip

We  propose  the converse of the preceding theorem as a conjecture.
\begin{centreconj}
If a simplex has a unique centre, then it is equifacetal.
\end{centreconj}
As noted earlier, this conjecture is true in dimension $3$.

Devid\'{e} \cite{Devide1975} proved in general that the incentre, centroid, and circumcentre of an equifacetal simplex must coincide, and asked whether the converse would hold in dimension $4$.  

A proof of the conjecture in higher dimensions would appear to require a larger family of centres, going beyond the classical ones. In forthcoming work the author will verify the conjecture for a suitably weak definition of ``centre'' (not necessarily defined for all simplices) and will verify the conjecture for a strong definition  of ``centre'' (using explicit centres defined for all simplices) in dimensions up to at least 6.

\section{The combinatorial type of a simplex}\label{sec:combinatorial}
Identify the $1$-skeleton of an $n$-simplex with $K_{n+1}$, the complete graph on $n+1$ vertices. We define the \emph{strong combinatorial type} of a general $n$-simplex to be the edge-coloured complete graph  $K_{n+1}$ on $n+1$ vertices, where the lengths of the edges are enumerated and a suitable label or colour is chosen for each edge of a given length.  In other words we forget the actual lengths of the edges and only remember which edges have the same length. We will observe that any colouring of $K_{n+1}$ arises as the combinatorial type of some $n$-simplex.
\subsection{Realization of given simplex edge lengths}
Consider the edge length function
\[
\mathcal{E}: (\mathbb{R}^{n})^{{n+1}}\to \mathbb{R}_{+}^{\binom{n+1}{2}}
\]
It is understood that the domain of $\mathcal{E}$ is those $(n+1)$-tuples $(A_{0}, A_{1}, \dots , A_{n})$ that actually form an $n$-simplex, which is given by a suitable determinant non-vanishing: 
\[
\det
\left(
\ 
\begin{matrix}
A_{1}-A_{0}&
A_{2}-A_{0}&
\dots&
A_{n}-A_{0}&
\end{matrix}
\right)
\ne 0
\]
It is an interesting question how to characterize the image of the edge length function.  When $n=1$ any edge length arises.  When $n=2$ a triple of possible edge lengths is realizable if and only if the requisite triangle inequalities are satisfied.   In general more conditions are required than just the triangle inequalities for all the faces, as one can easily see.  

This question was first answered by Menger \cite{Menger1931} (See the Third Fundamental Theorem on p. 738) in terms of inequalities for various determinants. A complete answer in the current form goes back to Schoenberg \cite{Schoenberg1935}, where it was shown that realizability is equivalent to a certain square matrix whose entries are polynomials in the proposed edge lengths being positive definite.   See also Blumenthal \cite{Blumenthal1953} for an extended treatment of related problems. A version of the result is also given in Dekster and Wilker \cite{DeksterWilker1987} and in Berger \cite{Berger1987}, Theorem 9.7.3.4. Rivin \cite{Rivin2003} has given a recent clear treatment in modern terms and in the present spirit.

Suppose given a set of $n(n+1)/2$ positive real numbers $\{  \ell_{ij}\}$, which are intended edge lengths of an $n$-simplex. Here we understand that $0\le i,j \le n$, $\ell_{ii}=0$, and $\ell_{ij}=\ell_{ji}$. Form the symmetric $n$ by $n$ matrix $M=[m_{ij}]$, where $m_{ij}=\ell_{i0}^{2}+\ell_{j0}^{2}-\ell_{ij}^{2}$.

\begin{thm}\label{posdef}
There is an $n$-simplex $A_{0}A_{1}\dots  A_{n}$ with $||A_{i}-A_{j}||=\ell_{ij}$ if and only if the matrix $M$ is positive definite. 
\qed
\end{thm}
See the references above for a proof. For present purposes we content ourselves with following immediate consequence, which is needed for what follows.
\begin{cor}
The edge length function has open image.
\qed
\end{cor}
In other words, given an $n$-simplex $S=A_{0}A_{1}\dots  A_{n}$, and given a set of positive real numbers $\{  \ell_{ij}\}$, the numbers are realized as the edge lengths of a simplex if the $\ell_{ij}$ are sufficiently close to the edge lengths $|A_{i}-A_{j}|$ of the given simplex.

We conclude this part with the observation that any edge-colouring of the complete graph $K_{n+1}$ arises as the combinatorial type of an $n$-simplex.
\begin{thm}\label{realize}
Suppose the edges of a complete graph on $n+1$ vertices are arbitrarily coloured.
Then, arbitrarily close to an equilateral $n$-simplex, there is an $n$-simplex whose edge graphs of constant length are precisely the subgraphs of given colour. in other words, every coloured graph decomposition of the complete graph on $n+1$ vertices arises from an $n$-simplex.
\end{thm}
\emph{Proof.  } 
For each colour choose a distinct real number arbitrarily close to $1$.  Then the openness of the space of realizable edge lengths shows that there is a simplex with these prescribed edge lengths close to the equilateral simplex with all edges of length $1$.
\qed\bigskip

Now we need a way of describing the strong combinatorial types that correspond  to equifacetal $n$-simplices.
\section{Restrictions on equifacetal simplices in higher dimensions}\label{sec:restrictions}
Here we observe the restrictions on the strong combinatorial type of an $n$-simplex that are imposed by being equifacetal.

As above we associate to each edge length $\ell$ of an equifacetal $n$-simplex $S$ the graph $G_{\ell}$ of all edges of that length. In this way we obtain the strong combinatorial type, or what we also call the \emph{coloured graph decomposition} $K_{n+1}=G_{\ell_{1}}\cup G_{\ell_{2}}\cup \dots \cup G_{\ell_{r}}$.  We forget the actual length of the edges and only keep track of which edges have the same length. 
From the present point of view it is the transitivity of the isometry group that provides important information. The  following three assertions are immediate consequences of the vertex transitivity of the isometry group but are useful for recognizing equifacetal or non-equifacetal simplices.
\begin{cor}
An equifacetal $n$-simplex contains at most $n$ different edge lengths.
\end{cor}
\emph{Proof.  } 
A vertex has degree $n$ and each edge length graph contains every vertex.
\qed\bigskip

\begin{cor}
A scalene $n$-simplex, $n\ge 3$, is not a face of an equifacetal $(n+1)$-simplex.
\end{cor}
\emph{Proof.  } 
Each face of an $(n+1)$-simplex has ${{n+1}\choose 2} =\frac{n(n+1)}{2}$ edges.  If $n\ge 3$, then ${{n+1}\choose 2} > n+1$.
\qed\bigskip

\begin{thm}[Complementarity condition]\label{complementarity}
If $S=A_0A_1\dots A_n$ is an equifacetal $n$-simplex with associated coloured graph decomposition
$$K_{n+1}=G_{\ell_{1}}\cup G_{\ell_{2}}\cup \dots \cup G_{\ell_{r}}$$ and $v$ and $w$ are two vertices of $K_{n+1}$, then the two coloured graphs obtained by deleting $v$ or $w$ (and all edges incident with them)
$$K_{n+1}-v=G_{\ell_{1}}-v\cup G_{\ell_{2}}-v\cup \dots \cup G_{\ell_{r}}-v$$ and
$$K_{n+1}-w=G_{\ell_{1}}-w\cup G_{\ell_{2}}-w\cup \dots \cup G_{\ell_{r}}-w$$
are isomorphic as coloured graphs.
\qed
\end{thm}
In words, the coloured graph decomposition of the complete graph associated with an equifacetal $n$-simplex has the property that the coloured graph obtained by deleting any vertex is uniquely defined up to isomorphism of coloured graphs. 
The following converse statement holds as  well, by Theorem \ref{realize}.
\begin{thm}
A coloured graph  decomposition of $K_{n+1}$ that  satisfies the complementarity condition arises as the strong combinatorial type of an equifacetal simplex.
\qed
\end{thm}
In dimension 2 there is only the equilateral triangle, which yields a graph with three nodes each of degree 2.  

In dimension 3 there are three possible graphs for equifacetal tetrahedra: complete graph on 4 vertices (equilateral faces); one circular chain and one disconnected graph  containing two edges (isosceles faces); or three disconnected graphs, each containing two edges (all three edges in a face of distinct length).  And all these possibilities  occur.
\section{The bookkeeping partition}\label{sec:partition}
One can use the following bookkeeping device to keep track of these graphs.  To an equifacetal $n$-simplex we can associate a partition $[n_1,n_2,\dots, n_k]$ of $n$, where we include $n_i$ in the list whenever there is one of the graphs of one edge length, of vertex degree $n_i$.   The bookkeeping partition is actually an invariant of what we might call \emph{vertex-uniform colourings,} in which all vertices have the same numbers of incident edges of each colour. In particular, as we will see, the complementarity condition is not perfectly reflected in the bookkeeping partition.

We will refer to this bookkeeping partition $[n_1,n_2,\dots, n_k]$ of $n$ as the \emph{weak combinatorial type}  of an equifacetal $n$-simplex,  or, more generally, of a vertex-uniform $n$-simplex.

In even dimensions there is  an additional  restriction on the partition.

\begin{lem}[Parity condition]If $S=A_0A_1\dots A_n$, $n$ \emph{even}, is an equifacetal $n$-simplex,  then the associated book-keeping partition has all even entries.
\end{lem}
  
\emph{Proof.  } 
Let $G_{\ell}$  be the graph of edges of  length $\ell$. We need to show that  all vertices have  even degree. Since $n$ is even, the number vertices, $n+1$, is odd.  If each vertex of $G_{\ell}$ has odd degree, the the  total number  of  edge  ends is odd, contradicting the obvious fact that a graph has an even number of edge ends.
\qed\bigskip

So in dimension 2 we only have $[2]$.  In dimension 3 we have $[3]$, $[2,1]$, or $[1,1,1]$.  In dimension $4$ we find $[4]$ or $[2,2]$.  In dimension 5 we have $[5]$, $[4,1]$, $[3,2]$, $[3,1,1]$, $[2,2,1]$, $[2,1,1,1]$, and $[1,1,1,1,1]$.  And so forth.  

In general, in higher dimensions, a given partition can be realized in more than one way, as we will see in Section \ref{sec:examples}.

Devid\'{e} \cite{Devide1975} completely  characterized the actual pairs of edge lengths that can be realized from the  essentially unique strong combinatorial type with two colours in dimension 4, corresponding  to  the partition $[2,2]$. From the present point of view, one just  needs to plug the  two edge lengths $a$ and $b$ into the matrix $M$ of Theorem \ref{posdef} and interpret the condition that the matrix be positive definite.
\section{Realizing the complementarity condition}\label{sec:realizing}
It is natural to conjecture that every such partition $\pi$ of the dimension $n$, satisfying the  parity condition, arises from an equifacetal simplex.  It turns out to be true when $n$ is even, but not always so when $n$ is odd.

Here we explore some fundamental realizability results and provide one basic  non-realizability result. We also formulate a general conjecture about realizability.
\begin{thm}\label{thm:evencase}
Let $\pi=[n_{1}, n_{2}, \dots, n_{k}]$ be an even partition of the even integer $n$.  Then $\pi$ is realized by a colouring of the complete graph $K_{n+1}$ in such a way as to satisfy the complementarity condition, and hence by an equifacetal $n$-simplex.
\end{thm}
\emph{Proof.  } 
Note that if the result is true for a partition of $k$ terms $[n_1,n_2,\dots, n_k]$, then it is true for the partition $[n_1,n_2,\dots, n_{k-1}+n_k]$, as one sees by simply colouring all edges previously coloured $k-1$ or $k$ by a single colour.  Thus it suffices to realize the partition $\pi=[2,2,\dots , 2]$ with $n/2$ terms.  To this end let the vertices be $\{ 0,1,2,\dots,n \}$.  For each $k$, $1\le k\le n/2$, define
$$X_{k}=\{ij:|i-j|\equiv k \mod n+1\}$$
Each edge $ij$ appears in exactly one $X_{k}$ and all the vertices appear in each $X_{k}$.  Each vertex of $X_{k}$ has degree $2$ in $X_{k}$ since $n+1$ is odd.  Thus $X_{k}$ can be expressed as a union of cycles, showing that $\pi$ can be realized.  The sets $X_{k}$  are each invariant under a cyclic permutation of the vertices, which is transitive.  The complementarity condition follows.
\qed\bigskip

We will see how  far one can push  the preceding  argument in  the case when $n$ is odd.

\begin{thm}Let $\pi=[n_{1}, n_{2}, \dots, n_{k}]$ be a partition of the odd integer $n$ such that only  one $n_{i}$ is odd. Then $\pi$ is realized by a colouring of the complete graph $K_{n+1}$ in such a way as to satisfy the complementarity condition, and hence by an equifacetal $n$-simplex.
\end{thm}
\emph{Proof.  } 
As before it suffices to realize the partition $\pi=[2,2,\dots , 2, 1]$ with $(n+1)/2$ terms.  To this end let the vertices be $\{ 0,1,2,\dots,n \}$.  For each $k$, $1\le k\le (n+1)/2$, define
$$X_{k}=\{ij:|i-j|\equiv k \mod n+1\}$$
Each edge $ij$ appears in exactly one $X_{k}$ and all the vertices appear in each $X_{k}$.  Each vertex of $X_{k}$, $k < (n+1)/2$, has degree $2$ in $X_{k}$ since $k\not\equiv -k \mod n+1$ then.  Thus such $X_{k}$ can be expressed as a union of cycles, as before.  It remains  to account for  $X_{(n+1)/2}$.  This actually consists of $(n+1)/2$ disjoint arcs, where each vertex has degree $1$.  This realizes $\pi$.  Moreover, once again all the sets $X_{k}$  are each invariant under the cyclic permutation $i\to i+1\mod n+1$ of the vertices, which is transitive.  (The colouring is also invariant under the action of the dihedral group of order $2(n+1)$.) The complementarity condition thus holds.
\qed\bigskip

One can, in fact, split any $X_{k}$ above that consists of cycles that alternate odd and even numbered vertices. This corresponds to passing from the full automorphism group $D_{2(n+1)}$ to its transitive subgroup $D_{n+1}$. Certain of the 2 entries in the corresponding partition are changed to $1,1$. Careful accounting of which $X_{k}$ can be split leads to the following statement.
\par
Define
\begin{equation*}
\mathcal{O}(n)=
\begin{cases}
\frac{n+1}{2}, &\text{ if  }\frac{n+1}{2} \text{ is odd}\\
\frac{n+3}{2}, &\text{ if  }\frac{n+1}{2} \text{ is even}
\end{cases}
\end{equation*}

\begin{cor}\label{dihedralrepresentability}
Let $\pi=[n_{1}, n_{2}, \dots, n_{r}]$ be any partition of the odd integer $n$. Suppose the number of odd entries in $\pi$ is no more than \(\mathcal{O}(n)
\). Then $\pi$ is realized by a colouring of the complete graph $K_{n+1}$ in such a way as to satisfy the complementarity condition, and hence by an equifacetal $n$-simplex.
\qed
\end{cor}
This can be easily used in particular to show that all partitions of $n=5$, except the partition $[1,1,1,1,1]$, which will be ruled out below, can be realized by a colouring of  the complete graph $K_{6}$ in such a way as to satisfy the complementarity condition.

Here is the simplest result showing that certain partitions are not realizable in the odd case.

\begin{thm}\label{notones}
If $n$ is odd and the partition $[1,1,\dots,1]$ is realizable by an equifacetal $n$-simplex, then $n+1$ is a power of $2$.
\end{thm} 
\emph{Proof.  } 
If this partition is realized by an equifacetal $n$-simplex, then consider the isometry group $G$, which acts transitively on the vertices. We claim that the isotropy group at each vertex must be trivial. Suppose to the contrary that the isotropy group $G_{v}$ of some vertex is nontrivial.  Let $w$ be a vertex for which there is a $g\in G_{v}$ such that $gw\ne w$. Then the edges $[v,w]$ and $[gv,gw]=[v,gw]$ have the same colour, but are adjacent, contradicting the form of the given partition.  Thus we may conclude that the vertex isotropy group is trivial and $G$ has order $|G|=n+1$ and acts freely on the set of vertices.

Now let $g\in G$ be a nontrivial element and suppose that $gv=w$. Then we can conclude that we also have $gw=v$. For otherwise $[v,w]$ and $[gv,gw]=[w,gw]$ gives a pair of adjacent edges with the same length.  It follows that $g^{2}$ is the identity, for all $g\in G$.

Elementary group theory then implies that $G\approx (\mathbf{Z}/2)^{r}$ for some positive integer $r$. Thus $n+1=|G|=2^{r}$, as required. 
\qed\bigskip

It would be good to close the gap between  the existence and non-existence results. We therefore pose the following.
\begin{partitionprob}
Let $n$ be a  positive integer. Determine which partitions of $n$  are realizable by an equifacetal $n$-simplex.
\end{partitionprob}
The  work of this paper solves the problem for the case when $n$ is even. For $n$ odd the preceding results solve it for $n\le 5$. We expect to return to this problem in a subsequent paper.

\begin{remark}
Interestingly enough it turns out when $n$ is odd that standard results from graph theory show that any partition arises from a vertex-uniform coloured graph decomposition of  $K_{n+1}$ (with the property that all vertices have the same number of incident edges of each length or colour). When $n$ is odd, then $K_{n+1}$ has a decomposition into $(n-1)/2$ Hamiltonian cycles (of length $n+1$) and a $1$-factor.  See \cite{HartsfieldRingel1994}, Theorem 2.3.2.  This realizes the partition $[2,2,\dots,2,1]$.  Since $n+1$ is even, the Hamiltonian cycles may be decomposed into pairs of $1$-factors to realize $[1,1\dots,1]$. The non-realizability result of Theorem \ref{notones} thus shows that such a graph colouring can be vertex-uniform without being vertex-transitive.

When $n$ is even, vertex uniformity (as well as the stronger condition of vertex transitivity) requires that the entries in the associated partition be even. Therefore we know that any partition of even parity arises in this case, and only such partitions arise as vertex uniform colored graph decompositions of $K_{n+1}$.
\end{remark}

\section{A purely group-theoretic perspective}\label{sec:group}
We observe that any finite group $G$ of order $m$ acts by isometries on an equifacetal $n$-simplex, where $n=m-1$. To see this, consider the regular representation $G\to S_{n+1}$ given by the action of $G$ on itself by left translation.  Then $G$ acts on the set of edges, i.e. unordered pairs $[g_{1},g_{2}]$, where $g_{i}\in G$, $g_{1}\ne g_{2}$, by $g[g_{1},g_{2}]=[gg_{1},gg_{2}]$. Color each orbit of edges with a distinct colour. Clearly $G$  acts  by isometries on any  equifacetal realization of this strong combinatorial type. The full isometry group is harder to pin down.  The corresponding partition has the form $[2,\dots, 2, 1,\dots, 1]$, where the $1$'s correspond to conjugacy classes of involutions in $G$. In other words, edge orbits either consist of a disjoint union of cycles (of equal length) or consist of a pairwise disjoint union of edges.

For example, if $G=C_{m}$, the cyclic group of order $m=n+1$, we get the standard construction we used to in our basic realizability results. And the full isometry group would be the dihedral group of order $2(n+1)$.

More generally start with any finite group $G$ and subgroup $H< G$ of index $m$, and let $n=m-1$.  Then the left action of $G$ on $G/H$ induces a representation $G\to S_{n+1}$ and similarly gives rise to an equifacetal simplex. This action has a kernel $K=\cap_{g\in G}\ gHg^{-1}$. The full automorphism group is an appropriate super-group of $G/K$.
\section{Examples and questions}\label{sec:examples}
Here we mention a number of examples in dimensions $5$ and $6$ that show the limitations in general constructions.

\begin{example}
When $n=5$ the vertex isotropy group $G_{v}$ may not be transitive on edges incident at vertex $v$ and of the same length. For this we consider a colouring of $K_{6}$ and a corresponding equifacetal $5$-simplex that  arise from the partition $[1,1,1,2]$ when $n=5$. Let the vertices be $0,1,\dots,5$. Let edges $01, 23, 45$ have colour $1$; edges $12, 34, 50$ have colour $2$; edges $03, 14,25$ have colour $3$; and $02,24,40; 13,35,51$ have colour $4$. Then clearly no isometry fixes $0$ and maps $2\to 4$. By employing a reflection that has no fixed vertices we can see that the isometry group is indeed transitive on edges.
\end{example}
\begin{example}
When $n=5$ there need not be an isometry that interchanges two given vertices. For, in the preceding example, the edges of the triangles $02,24,40$; and $13,35,51$ cannot be reversed.
\end{example}
\begin{example}
When $n=5$ there are inequivalent realizations of the partition $[2,3]$ by equifacetal simplices and coloured graph decompositions. Here is  a realization in which the edges coloured $2$ (corresponding to the partition term $2$) form a connected cycle: $01,12,23,34,45,50$, in which those coloured $3$  are all the other  edges.  On the other hand one could equally well colour  the  edges $02,24,40$ and $13,35,51$ ``2'', and the remaining edges ``3''. A simple  drawing reveals that both of these coloured  graph decompositions satisfy the complementarity condition  and arise from equifacetal $5$-simplices. But  they are manifestly inequivalent.
\end{example}
\begin{example}
When $n=5$ there are vertex-uniform coloured graph decompositions realizing the partition $[1,2,2]$ that do not satisfy the complementarity condition (as well as realizations that do satisfy it). Consider  the following colouring: colour 1 is given by $04,13,25$; colour $2$ is given by $02,24,41,15,53,30$; and colour $3$ is given by $01,12,23,34,45,50$. Again a simple sketch figure reveals the non-complementarity, in which the end  points of the  colour $3$ arc in $K-0$ and $K-5$, respectively, are connected by  an edge with colour 1 in the second case, but not in the first case. It is easy to apply Theorem \ref{bestrealizability} to construct a colouring with partition $[2,3]$ that does satisfy complementarity, as well.
\end{example}
Based on our low-dimensional constructions it is plausible to suppose that the isometry group of an equifacetal $n$-simplex is transitive on the sets of edges of given length.  This in fact true for dimensions up through $5$.  General thinking shows that the set of edges of given length is certainly a union of orbits of the action of the isometry group on the set of edges. 

\begin{example}
When $n=6$ the standard realization of the partition $[4,2]$ has isometry group $D_{14}$ which does not act transitively on the set of  ``4'' edges. The standard realization of $[2,2,2]$ has the same isometry group.  When two of the sets of different colours are repainted to have the same colour, the isometry group, in this case, does not increase in compensation. See the appendix for a few more details.
\end{example}
In the introduction we observed that any segment is the edge of an equilateral triangle (Euclid) and a triangle is a facet of an equifacetal tetrahedron if and only if the triangle is acute. 

Is there an effective and straight-forward criterion for when an $(n-1)$-simplex is a facet of an equifacetal $n$-simplex in general?

Devid\'{e} \cite{Devide1975} gave such a classification in dimension $4$, aided by the fact that there is an essentially unique combinatorial type for a non-equilateral, equifacetal $4$-simplex. Unfortunately such a result seems to be elusive in higher dimensions, partly because of the example of the  partition $[1,1,1,1,1]$ that can be realized by a edge-colouring of $K_{6}$, so that all vertices look the same, with all $5$ different colours incident to every vertex, but such that the complementarity condition fails. Similarly we saw that there are colourings associated with the partition $[1,2,2]$ that do not satisfy the complementarity property, even though there are other colourings that do satisfy it.

One can fairly easily develop a result that characterizes colourings of $K_{n}$ that extend to vertex-uniform colourings of $K_{n+1}$. Such an extension is essentially unique. Then one can inspect the resulting extension to see if it also satisfies the complementarity property (vertex transitivity). Surely there is a better result available.
In relatively low dimensions we can list the combinatorial types of equifacetal simplices and then delete a vertex.  But that process naturally becomes unwieldly as the dimension increases.

\section{Combinatorial types in low dimensions}\label{sec:table}
In each dimension $n$ there is of course the trivial partition $[n]$ and its unique corresponding coloured graph, the complete graph $K_{n+1}$, with all edges coloured the same, which is realized by a standard equilateral $n$-simplex, with isometry group $S_{n+1}$. We enumerate the possible nontrivial partitions and their possible coloured graph realizations, together with the corresponding isometry groups.

In small dimensions one can often deduce the full isometry group as follows. If there is an edge orbit of type $2$ that is a single cycle of $n+1$ edges, then the isometry group overall is a subgroup of the isometry group of the cycle. On the other hand the isometry group must be at least of order $n+1$.  So in such a case the isometry group is cyclic $C_{n+1}$ of order $n+1$, dihedral $D_{2(n+1)}$ of order $2(n+1)$, or dihedral $D_{n+1}$ of  order $n+1$. If two distinct $1$-factors, each of type 1, fit together to form a single cycle, then the full isometry group must be just $D_{n+1}$.

\paragraph{Dimension $3$} There are two such partitions
\begin{enumerate}
\item $[1,1,1]$. Realized by a simplex in which the three pairs of disjoint  edges get different colours. The group is the dihedral group of order $4$, $D_{4}\approx C_{2}\times C_{2}$.
\item $[2,1]$. This is a specialization of the preceding case. It can be viewed as a square with its four edges coloured $a$, with two diagonals coloured $b$. The group is $D_{8}$. 
\end{enumerate}

\paragraph{Dimension $4$} There is one such partition satisfying the parity condition.
\begin{enumerate}
\item $[2,2]$. Either type $2$ term must correspond to a $5$-cycle passing through all the vertices. Think of one of them as being realized by a pentagon coloured $a$. Invariance under a cyclic or dihedral action shows that there is a unique way to add the second colour, as a pentagram coloured $b$. Compare Devid\'{e} \cite{Devide1975}. The group is the dihedral group of order $10$, $D_{10}$. 
\end{enumerate}

\paragraph{Dimension $5$} There are five such partitions. ( $[1,1,1,1,1]$ cannot be realized.) Note that the realizations of $[2,1,1,1]$ and of $[3,1,1]$ both have the same group, namely $D_{6}$.
\begin{enumerate}
\item $[2,1,1,1]$. Two sets of  type $1$ edges together must form a hexagon. So we have a hexagon with edges alternately coloured $a$ and $b$. The type $2$  edges cannot form a hexagon, for otherwise it could be split to yield a realization of $[1,1,1,1,1]$. We conclude that the type $2$ edges form a pair of triangles, each connecting alternate vertices of the hexagon. Color them $c$. Connect opposite vertices on the hexagon with edges coloured $d$ to realize the final type $1$ family. The group is $D_{6}$ coming from the symmetries of the $ab$ hexagon.
\item $[3,1,1]$. The two sets of type $1$ edges must form a hexagon. So we have a hexagon with edges alternately coloured $a$ and $b$, with the remaining edges coloured $c$. The group is $D_{6}$ coming from the symmetries of the $ab$ hexagon.
\item $[2,2,1]$. One set of type $2$ edges forms a hexagon, coloured $a$ and the other set of type $2$ edges form two disjoint triangles, each connecting alternate vertices of the hexagon. (Invariance under a vertex transitive group shows that they cannot both be hexagons or both be pairs of triangles.) Opposite vertices on the hexagon are connected by edges coloured $c$.  The group is $D_{12}$.
\item $[3,2]$. Here there are two coloured graph configurations, distinguished by their groups, one of order $24$ and the other of order $72$. They can be viewed as the two possible specializations of the realization of $[2,2,1]$ above, in which the colour $c$ is changed to $a$ or to $b$.
\begin{enumerate}
\item A hexagon with its $6$ edges coloured $a$ and the remaining $9$ edges, connecting all of  the other possible pairs of vertices, coloured $b$. The group is $D_{12}$.
\item Two disjoint triangles with edges coloured $a$ and the remaining $9$ edges, connecting all of the other possible pairs of vertices, coloured $b$. The group is $D_{6}\times D_{6}=S_{3}\times S_{3}$ extended by $C_{2}$. The two triangles can be independently mapped into themselves, or they can be interchanged.
\end{enumerate}

\item $[4,1]$. Three disjoint edges are coloured $a$ and all the other edges of $K_{6}$ are coloured $b$. The group is $C_{2}\times C_{2}\times C_{2}$, extended by $S_{3}$, since the three $a$-segments must be preserved but can be independently reversed and also permuted among themselves.
\end{enumerate}

\paragraph{Dimension $6$} There are two such partitions satisfying the parity condition. They both have the same group, $D_{14}$.
\begin{enumerate}
\item $[4,2]$. The type $2$ edges must form a single $7$-gon coloured $a$, say, since they form a union of cycles all of the same size. The remaining edges of $K_{7}$ are coloured $b$. The group is $D_{14}$. Comparing with the realization of $[2,2,2]$  below it becomes apparent that $D_{14}$ does not act transitively on the set of $14$ type $4$ edges, but rather has two orbits of size $7$ there.
\item $[2,2,2]$. As above we have three $7$-gons. The first one may be coloured $a$ with vertices labelled $0$ through $7$. We must have a vertex transitive group that preserves this $7$-gon, so it contains $C_{7}$, for which a generator can be viewed as mapping $i\to i+1\mod 7$. The other two $7$-gons are then $0246135$ (skip one around) and $0362514$ (skip 2 around). The full group of isometries is clearly $D_{14}$. In this case $D_{14}$ does act transitively on each set of $7$ edges of each colour.
\end{enumerate}

\noindent
Professor Allan L. Edmonds,\\
Department of Mathematics,\\
Indiana University,\\
Bloomington, IN 47405 USA.\\

E-mail: edmonds@indiana.edu\\

52B11: CONVEX AND DISCRETE GEOMETRY; $n$-dimensional polytopes.

 \end{document}